\documentclass{elsart}
\usepackage{amssymb}
\usepackage{graphicx}

\begin{document}

\begin{frontmatter}

\title{Limit~Cycles~of~a~Quadratic~System~with
Two~Parallel Straight~Line-Isoclines\thanksref{label1}}
\thanks[label1]{This work is supported by the Netherlands Organization for Scientific
Research (NWO).}

\author{Valery A. Gaiko\thanksref{label2}}
\ead{vlrgk@yahoo.com}
\thanks[label2]{The author is grateful to the Delft Institute of Applied Mathematics 
of TU~Delft for hospitality during his stay at the University in the period of October, 
2007~-- March, 2008. He is also thankful to Wim van Horssen and John Reyn
for very productive discussions.}

\address{Belarusian State University of Informatics and Radioelectronics,
Department~of~Mathematics,~L.\,Beda~Str.\,6--4,~Minsk~220040,~Belarus}

\begin{abstract}
In this paper, a quadratic system with two parallel straight line-isoclines is considered.
This system corresponds to the system of class~II in the classification of Ye Yanqian \cite{Ye}.
Using the field rotation parameters of the constructed canonical system and geometric 
properties of the spirals filling the interior and exterior domains of its limit cycles,
we prove that the maximum number of limit cycles in a~quadratic system with two 
parallel straight line-isoclines and two finite singular points is equal to two. Besides, 
we obtain the same result in a different way: applying the Wintner--Perko termination 
principle for multiple limit cycles and using the methods of global bifurcation theory 
developed in \cite{Gaiko}.

\end{abstract}

\end{frontmatter}

\section{Introduction}
\label{1}

We consider the system of differential equations
    $$
    \begin{array}{l}
\dot{x}=a_{00}+a_{10}x+a_{01}y+a_{20}x^{2}+a_{11}xy+a_{02}y^{2},\\ \dot{y}=b_{00}+b_{10}x+b_{01}y+b_{20}x^{2}+b_{11}xy+b_{02}y^{2}\\
    \end{array}
    \eqno(1.1)
    $$
with the real coefficients $a_{ij},$ $b_{ij}$ in the real variables $x,\:y,$ where
at least one quadratic term has a coefficient unequal to zero. Such a system will
be referred to as a quadratic system. The main problem of the qualitative theory 
of system (1.1) is \emph{Hilbert's Sixteenth Problem} on the maximum number 
and relative position of its limit cycles, i.\,e., closed isolated trajectories of (1.1) 
\cite{BL}, \cite{Gaiko}, \cite{Ilyash}--\cite{Ye}. The solution of this problem 
could give us all possible phase pertraits of system (1.1) and, thus, could complete 
its qualitative analysis.
   \par
Earlier \cite{gai1}--\cite{gai5}, we studied a quadratic system with two intersecting 
straight line-isoclines. In this paper, we will study a quadratic system with two parallel 
straight line-isoclines. Such a system corresponds to the system of class~II in the 
classification of Ye Yanqian~\cite{Ye}:
$$
    \begin{array}{l}
\dot{x}=\lambda\,x-y+l\,x^{2}+m\,xy+n\,y^{2},\\ 
\dot{y}=x+x^{2}.
    \end{array}
    \eqno(1.2)
    $$
Applying a new geometric approach to the study of limit cycle bifurcations developed 
in~\cite{gai5}, we will prove that a quadratic system with two parallel sraight line-isoclines 
and two finite singular points has at most two limit cycles (a~quadratic system with one 
sraight line-isocline was studied in detail in \cite{BL}). The same result will be obtained 
in a different way: applying the Wintner--Perko termination principle for multiple limit 
cycles~\cite{Perko} and using the methods of global bifurcation theory developed 
in~\cite{Gaiko}.
   \par
In particular, in Section~2, we construct a canonical system with field rotation 
parameters corresponding to the system of class~II in the classification of 
Ye Yanqian \cite{Ye}. In Section~3, using the canonical system and geometric 
properties of the spirals filling the interior and exterior domains of limit cycles, 
we obtain the main result of this paper on the maximum number of limit cycles
of a quadratic system with two parallel straight line-isoclines. In Section~4, 
we obtain the same result applying the Wintner--Perko termination principle 
for multiple limit cycles. 

\section{Canonical systems}
\label{2}

First, we have to construct canonical quadratic systems with field rotation 
parameters for studying limit cycle bifurcations. The following theorem is valid.
    \par
    \medskip
\noindent\textbf{Theorem 2.1.}
    \emph{A quadratic system with limit cycles can be reduced to
the canonical form}
    $$
\begin{array}{l}
    \dot{x}=-y\,(1+x+\alpha\,y),\\
    \dot{y}=x+(\lambda+\beta+\gamma)y+a\,x^{2}
        +(\alpha+\beta+\gamma)xy+c\,\gamma\,y^{2}
        \end{array}
    \eqno(2.1)
    $$
\emph{or}
    \\[-4mm]
    $$
\begin{array}{l}
    \dot{x}=-y\,(1+\nu\,y)\equiv P, \quad \nu=0;1,\\
    \dot{y}=x+(\lambda+\beta+\gamma)y+a\,x^{2}
        +(\beta+\gamma)xy+c\,\gamma\,y^{2}\equiv Q.
\end{array}
    \eqno(2.2)
    $$
    \par
\noindent\textbf{Proof.} As was shown in \cite{Gaiko}, by means of Erugin's
two-isocline method~\cite{erug}, an arbitrary quadratic system with limit cycles
can be reduced to the form
    $$
\begin{array}{l}
    \dot{x}=-y+mxy+ny^{2},\\
    \dot{y}=x+\lambda y+ax^{2}+bxy+cy^{2},
\end{array}
    \eqno(2.3)
    $$
where $m=-1$ or $m=0.$
    \par
Input the field rotation parameters into this system so that (2.1)
corresponds to the case of $m=-1$ and (2.2) corresponds to the
case of $m=0.$
    \par
Compare (2.1) with (2.3) when $m=-1.$ Firstly, we have changed
several parameters: $n$ by $-\alpha;$ $b$ by $\beta;$ $c$ by
$c\,\gamma.$ Secondly, we have input additional terms into the
expression for $\dot{y}\!:$ $(\beta+\gamma)\,y$ and
$(\alpha+\gamma)\,xy.$ Similar transformations have been made in
system (2.3) when $m=0;$ but in this case, we have denoted $n$ by
$\nu$ assigning two principal values to this parameter: 0 and 1. 
It is obvious that all these transformations do not restrict generality 
of systems (2.1) and (2.2) in comparison with system (2.3), which 
proves the theorem.\qquad $\Box$
    \medskip
    \par
System (2.1) is a basic system for studying limit cycle
bifurcations. It contains four field rotation parameters:
$\lambda,\:\alpha,\:\beta,\:\gamma.$ This system has been 
considered in \cite{gai5}. Now we will consider system (2.2).
The following lemma is valid for (2.2).
    \par
    \medskip
\noindent\textbf{Lemma 2.1.}
    \emph{Each of the parameters $\lambda,\:\beta,$ and $\gamma$
rotates the vector field of system $(2.2)$ in the domains of existence 
of its limit cycles, under the fixed other parameters of this system,
namely$\,:$ when the parameter $\lambda,\:\beta,$ or
$\gamma$ increases $($decreases$),$ the field is rotated in
positive $($negative$)$ direction, i.\,e., counterclockwise
$($clockwise$),$ in the do\-mains, respectively$\,:$}\\[-2mm]
    $$
    1+\nu\,y<0 \ (>0);
    $$
    $$
    (1+x)(1+\nu\,y)<0 \ (>0);
    $$
    $$
    (1+x+c\,y)(1+\nu\,y)<0 \ (>0).
    $$
    \par
\noindent\textbf{Proof.} \ Using the definition of a field
rotation parameter~\cite{duff},~\cite{Gaiko}, we can calculate the
following determinants:\\
    $$
\Delta_{\lambda}=PQ'_{\lambda}-QP'_{\lambda}=-y^{2}(1+\nu\,y);
    $$
    $$
\Delta_{\beta}=PQ'_{\beta}-QP'_{\beta}=-y^{2}(1+x)(1+\nu\,y);
    $$
    $$
\Delta_{\gamma}=PQ'_{\gamma}-QP'_{\gamma}=-y^{2}(1+x+c\,y)(1+\nu\,y).
    $$
    \par
Since, by definition, the vector field is rotated in positive direction 
(counter\-clock\-wise) when the determinant is positive and in 
negative direction (clock\-wise) when the determinant is 
negative~\cite{duff},~\cite{Gaiko} and since the obtained domains
correspond to the domains of existence of limit cycles of (2.2),
the lemma is proved.\qquad $\Box$

\section{Limit cycle bifurcations}
\label{3}

We will study limit cycle bifurcations of canonical system (2.2) with two parallel
straight line-isoclines and three field rotation parameters. We will consider the
case when system (2.2) has only two finite singularities: a saddle and an anti-saddle
(all other cases can be considered in an absolutely similar way). Let us prove the 
following theorem.
    \par
    \medskip
\noindent\textbf{Theorem 3.1.}
    \emph{System (2.2) with two parallel straight line-isoclines and two finite 
singular points can have at least two limit cycles surrounding the origin.}
    \medskip
    \par
    \noindent\textbf{Proof.} \
To prove the theorem, fix, for example, $a=1$ and take $c>1$ in system (2.2) for $\nu=1.$
Then vanish all field rotation parameters of (2.2), $\beta=\gamma=\lambda=0\!:$ 
    $$
\begin{array}{l}
    \dot{x}=-y\,(1+y),\\
    \dot{y}=x+x^{2}.
\end{array}
    \eqno(3.1)
    $$
We have got a system with the zero divergence and four finite singular points: 
two centers and two saddles (a Hamiltonian case).
    \par
Input, for example, a positive parameter $\gamma$ into system (3.1):
    $$
\begin{array}{l}
    \dot{x}=-y\,(1+y),\\
    \dot{y}=x+\gamma\,y+x^{2}+\gamma\,xy+c\,\gamma\,y^{2}.
\end{array}
    \eqno(3.2)
    $$
On increasing the parameter $\gamma,$ the vector field of~(3.2) is rotated 
in negative direction (clockwise) and the center at the origin turns into an unstable focus.
   \par
Suppose that $\gamma$ satisfies the condition 
  $$
\displaystyle-1+2\left(c-\sqrt{c(c-1)}\right)<\gamma<-1+2\left(c+\sqrt{c(c-1)}\right).
  \eqno(3.3)
  $$ 
In this case we will have only two finite sigularities: 
a saddle $(-1,0)$ and an anti-saddle $(0,0).$ 
    \par
Fix $\gamma$ and input a negative parameter $\beta$ into system (3.2):
    $$
\begin{array}{l}
    \dot{x}=-y\,(1+y),\\
    \dot{y}=x+(\beta+\gamma)\,y+x^{2}+(\beta+\gamma)\,xy+c\,\gamma\,y^{2}.
\end{array}
    \eqno(3.4)
    $$
On decreasing the parameter $\beta,$ the vector field of~(3.4) is rotated in
positive direction, and, for some value $\beta_{S}$ of the parameter $\beta,$
a separatrix loop is formed around the origin generating a stable limit cycle for
$\beta<\beta_{S}$ (an~unstable limit cycle cannot appear from the origin
because of the negative first focus quantity at the origin for $\gamma>0$ 
when $\beta+\gamma=0$ \cite{BL}).
    \par
Fix $\beta$ satisfying the condition $0<\beta+\gamma\ll1$ and input 
a positive pa\-ra\-me\-ter $\lambda$ into system (3.4):
    $$
\begin{array}{l}
    \dot{x}=-y\,(1+y),\\
    \dot{y}=x+(\lambda+\beta+\gamma)\,y+x^{2}+(\beta+\gamma)\,xy+c\,\gamma\,y^{2}.
\end{array}
    \eqno(3.5)
    $$
To have still two finite singularities, we also suppose that the parameters $\beta,$~$\gamma,$
and~$\lambda$ satisfy the condition 
  $$
-1-\sqrt{4c\gamma-1}<\beta+\gamma+\lambda<-1+\sqrt{4c\gamma-1}.
  \eqno(3.6)
  $$ 
On increasing the parameter $\lambda,$ the vector field of~(3.5) is rotated in
negative direction, and, for some value $\lambda=\lambda_{S},$ a separatrix loop 
is formed around the origin again generating an unstable limit cycle for
$\lambda>\lambda_{S}$ (the stable limit cycle cannot disappear through the loop 
because of the positive divergence at the saddle $(-1,0)$ for $\lambda>0$ \cite{BL}).
   \par
Thus, we have obtained at least two limit cycles surrounding the focus $(0,0),$ 
which proves the theorem.\qquad $\Box$
    \par
Let us prove now a much stronger theorem (it is the main result of our paper).
    \par
    \medskip
\noindent\textbf{Theorem 3.2.}
    \emph{System (2.2) with two parallel straight line-isoclines and two finite 
singular points has at most two limit cycles surrounding the origin.}
    \medskip
    \par
\noindent\textbf{Proof.} \ Consider again system (2.2) for $\nu=1,$ $a=1,$ and $c>1$ 
supposing that condition (3.6) is also valid. All other particular cases of (2.2) can be 
considered in a similar way.
   \par
Vanishing all field rotation parameters of system (3.1), $\beta=\gamma=\lambda=0,$
we have got again Hamiltonian system (3.1) with four finite singular points: 
two centers and two saddles. Let us input successively the field rotation parameters into (3.1).
    \par
Begin, for example, with the parameter~$\gamma$ supposing that $\gamma>0.$ Then we
get system (3.2). On increasing the parameter $\gamma,$ the vector field of~(3.2) is rotated 
in negative direction (clockwise) and the center at the origin turns into an unstable focus. We 
also suppose that $\gamma$ satisfies condition (3.3) and that we have only two finite sigularities
in this case: a focus $(0,0)$ and a saddle $(-1,0).$ 
    \par
Fix $\gamma$ and input a negative parameter $\beta$ getting system (3.4).
On decreasing the parameter $\beta,$ the vector field of~(3.4) is rotated in
positive direction, and, for some value $\beta_{S}$ of the parameter $\beta,$
a separatrix loop is formed around the origin generating a stable limit cycle for
$\beta<\beta_{S}.$ As we noted above, a~limit cycle cannot appear from the 
origin because of the negative first focus quantity at the origin for $\gamma>0$ 
when $\beta+\gamma=0$ \cite{BL}.
    \par
Denote the limit cycle by $\Gamma\!_{1},$ the domain inside the
cycle by $D_{1},$ the domain outside the cycle by $D_{2}$ and
consider logical possibilities of the appearance of other 
(semi-stable) limit cycles from a ``trajectory concentration''
surrounding the focus $(0,0).$ It is clear that on decreasing
$\beta,$ a semi-stable limit cycle cannot appear in the domain
$D_{2},$ since the outside spirals winding onto the cycle will 
untwist and the distance between their coils will increase 
because of the vector field rotation in positive direction.
    \par
By contradiction, we can also prove that a semi-stable limit cycle
cannot appear in the domain $D_{1}.$ Suppose it appears in this
domain for some values of the parameters $\gamma^{*}>0$ and
$\beta^{*}<0.$ Return to initial system (3.1) and change the
order of inputting the field rotation parameters. Input first the
parameter $\beta<0\!:$
    $$
\begin{array}{l}
    \dot{x}=-y\,(1+y),\\
    \dot{y}=x+\beta\,y+x^{2}+\beta\,xy.
\end{array}
    \eqno(3.7)
    $$
Fix it under $\beta=\beta^{*}.$ The vector field of~(3.7) is rotated 
in negative direction and $(0,0)$ becomes a stable focus. Inputting 
the parameter $\gamma>0$ into (3.7), we have got again system (3.4), 
the vector field of which is rotated in positive direction. Under this rotation, 
for $\gamma=-\beta,$ the focus $(0,0)$ changes the character of its
stability, and a stable limit cycle $\Gamma\!_{1}$ appears from the origin.
This cycle will expand, the focus spirals will untwist, and the distance
between their coils will increase on increasing the parameter
$\gamma$ to the value $\gamma=\gamma^{*}.$ It follows that there
are no values of $\gamma=\gamma^{*}$ and $\beta=\beta^{*},$
for which a semi-stable limit cycle could appear in the domain
$D_{1}.$
    \par
This contradiction proves the uniqueness of a limit cycle
surrounding the focus $(0,0)$ in system (3.4) for any values of
the parameters $\gamma$ and $\lambda$ of different signs.
Obviously, if these parameters have the same sign, system (3.4)
has no limit cycles surrounding $(0,0)$ at all.
    \par
Let system (3.4) have the unique limit cycle $\Gamma\!_{1}.$ 
Fix the parameters $\gamma>0,$\linebreak $\beta<0$ and input the third
parameter, $\lambda>0,$ getting system (3.5). On increasing the 
parameter $\lambda,$ the vector field of~(3.5) is rotated in
negative direction, and, for some value $\lambda=\lambda_{S},$ 
a separatrix loop is formed around the origin again generating 
an unstable limit cycle for $\lambda>\lambda_{S}.$ Note again 
that a limit cycle cannot disappear through the loop because of 
the positive divergence at the saddle $(-1,0)$ for $\lambda>0$ 
\cite{BL}. Denote this cycle by $\Gamma_{2}.$ On further increasing 
$\lambda,$ the limit cycle $\Gamma_{2}$ will join with $\Gamma\!_{1}$ 
forming a semi-stable limit cycle, $\Gamma_{\!12},$ which will disappear 
in a ``trajectory concentration'' surrounding the origin $(0,0).$ Can another 
semi-stable limit cycle appear around the origin in addition to $\Gamma_{\!12}?$ 
It is clear that such a limit cycle cannot appear either in the domain $D_{1}$ 
bounded by the origin and $\Gamma\!_{1}$ or in the domain $D_{3}$ bounded 
on the inside by $\Gamma_{2}$ because of the increasing distance between the
spiral coils filling these domains on increasing $\lambda.$
    \par
To prove impossibility of the appearance of a semi-stable limit
cycle in the domain $D_{2}$ bounded by the cycles $\Gamma\!_{1}$
and $\Gamma_{2}$ (before their joining), suppose the contrary,
i.\,e., for some set of values of the parameters $\gamma^{*}>0,$
$\beta^{*}<0,$ and $\lambda^{*}>0,$ such a semi-stable cycle
exists. Return to system (3.1) again and input first the parameters
$\gamma>0$ and $\lambda>0\!:$
    $$
\begin{array}{l}
    \dot{x}=-y\,(1+y),\\
    \dot{y}=x+(\lambda+\gamma)\,y+x^{2}+\gamma\,xy+c\,\gamma\,y^{2}.
\end{array}
    \eqno(3.8)
    $$
Both parameters act in a similar way: they rotate the vector field of (3.7) 
in negative direction turning the origin $(0,0)$ into an unstable focus.
    \par
Fix these parameters under $\gamma=\gamma^{*},$
$\lambda=\lambda^{*}$ and input the parameter $\beta<0$ into
(3.8) getting again system (3.5). Since, in our assumption, this
system has two limit cycles for $\beta<\beta^{*},$ there exists
some value of the parameter, $\beta_{12}$
$(\beta^{*}<\beta_{12}<0),$ for which a semi-stable limit cycle,
$\Gamma_{12},$ appears in system (3.5) and then splits into a
stable cycle, $\Gamma\!_{1},$ and an unstable cycle, $\Gamma_{2},$
on further decreasing $\beta.$ The formed domain $D_{2}$
bounded by the limit cycles $\Gamma\!_{1},$ $\Gamma_{2}$ and
filled by the spirals will enlarge, since, by the properties of a
field rotation parameter, the interior stable limit cycle
$\Gamma\!_{1}$ will contract and the exterior unstable limit cycle
$\Gamma_{2}$ will expand on decreasing $\beta.$ The distance
between the spirals of the domain $D_{2}$ will naturally increase,
which will prohibit the appearance of a semi-stable limit cycle in
this domain for $\beta<\beta_{12}.$ Thus, there are no such
values of the parameters, $\gamma^{*}>0,$ $\lambda^{*}>0,$
$\beta^{*}<0,$ for which system (3.5) would have an additional
semi-stable limit cycle.
    \par
Obviously, there are no other values of the parameters $\lambda,$ 
$\beta,$ $\gamma,$ for which system (3.5) would have more than two 
limit cycles surrounding the origin $(0,0).$ It follows that system (3.5) 
and, hence, system (2.2) can have at most two limit cycles. 
The theorem is proved.
\qquad $\Box$

\section{The Wintner--Perko termination principle}
\label{4}

In \cite{Gaiko}, for the global analysis of limit cycle bifurcations, 
we used the Wintner--Perko termination principle
which was stated for relatively prime, pla\-nar, analytic systems
and which connected the main bifurcations of limit
cycles~\cite{Perko}. Let us formulate this principle
for the polynomial system
    $$
    \mbox{\boldmath$\dot{x}$}=\mbox{\boldmath$f$}
    (\mbox{\boldmath$x$},\mbox{\boldmath$\mu$)},
    \eqno(4.1_{\mbox{\boldmath$\mu$}})
    $$
where $\mbox{\boldmath$x$}\in\textbf{R}^2;$ \
$\mbox{\boldmath$\mu$}\in\textbf{R}^n;$ \
$\mbox{\boldmath$f$}\in\textbf{R}^2$ \ $(\mbox{\boldmath$f$}$ is a
polynomial vector function).
    \par
    \medskip
\noindent\textbf{Theorem 4.1 (Wintner--Perko termination principle).}
    \emph{Any one-para\-me\-ter fa\-mi\-ly of multiplicity-$m$
limit cycles of relatively prime polynomial system $(4.1_{\mbox{\boldmath$\mu$}})$ 
can be extended in a unique way to a maximal one-parameter family of multiplicity-$m$ 
limit cycles of $(4.1_{\mbox{\boldmath$\mu$}})$ which is either open or cyclic.}
    \par
\emph{If it is open, then it terminates either as the parameter or
the limit cycles become unbounded; or, the family terminates
either at a singular point of $(4.1_{\mbox{\boldmath$\mu$}}),$
which is typically a fine focus of multiplicity~$m,$ or on a
(compound,) separatrix cycle of $(4.1_{\mbox{\boldmath$\mu$}}),$ 
which is also typically of multiplicity~$m.$}
    \medskip
    \par
The proof of the Wintner--Perko termination principle for general
polynomial system $(4.1_{\mbox{\boldmath$\mu$}})$ with a vector
parameter $\mbox{\boldmath$\mu$}\in\textbf{R}^n$ parallels the
proof of the pla\-nar termination principle for the system
    $$
    \vspace{1mm}
    \dot{x}=P(x,y,\lambda),
        \quad
    \dot{y}=Q(x,y,\lambda)
    \eqno(4.1_{\lambda})
    \vspace{2mm}
    $$
with a single parameter $\lambda\in\textbf{R}$ (see \cite{Gaiko},
\cite{Perko}), since there is no loss of generality in assuming
that system $(4.1_{\mbox{\boldmath$\mu$}})$ is parameterized by a
single parameter $\lambda;$ i.\,e., we can assume that there
exists an analytic mapping $\mbox{\boldmath$\mu$}(\lambda)$ of
$\textbf{R}$ into $\textbf{R}^n$ such that
$(4.1_{\mbox{\boldmath$\mu$}})$ can be written as
$(4.1\,_{\mbox{\boldmath$\mu$}(\lambda)})$ or even
$(4.1_{\lambda})$ and then we can repeat everything, which had been
done for system $(4.1_{\lambda})$ in~\cite{Perko}. In particular,
if $\lambda$ is a field rotation parameter of $(4.1_{\lambda}),$
the following Perko's theorem on monotonic families of limit cycles
is valid.
    \par
    \medskip
\noindent\textbf{Theorem 4.2.}
    \emph{If $L_{0}$ is a nonsingular multiple limit cycle of
$(4.1_{0}),$ then $L_{0}$ belongs to a one-parameter family of
limit cycles of $(4.1_{\lambda});$ furthermore:}
    \par
1)~\emph{if the multiplicity of $L_{0}$ is odd, then the family
either expands or contracts mo\-no\-to\-ni\-cal\-ly as $\lambda$
increases through $\lambda_{0};$}
    \par
2)~\emph{if the multiplicity of $L_{0}$ is even, then $L_{0}$
bi\-fur\-cates into a stable and an unstable limit cycle as
$\lambda$ varies from $\lambda_{0}$ in one sense and $L_{0}$
dis\-ap\-pears as $\lambda$ varies from $\lambda_{0}$ in the
opposite sense; i.\,e., there is a fold bifurcation at
$\lambda_{0}.$}
    \par
    \medskip
In \cite{gai1}--\cite{gai5}, using Theorems~4.1 and~4.2, 
we have proved the following theorem.
    \par
    \medskip
\noindent\textbf{Theorem 4.3.}
    \emph{There exists no quadratic system having a swallow-tail
bifurcation surface of multiplicity-four limit cycles in its
pa\-ra\-meter space. In other words, a quadratic system cannot
have either a multi\-plicity-four limit cycle or four limit cycles
around a singular point (focus), and the maximum
multi\-plicity or the maximum number of limit cycles surrounding a
focus is equal to three.}
    \medskip
    \par
Applying the same approach, let us give an alternative proof
of Theorem~3.2.
   \par
\noindent\textbf{Proof (an alternative proof of Theorem~3.2).} \ 
The proof of this theorem is carried out by contradiction.
Consider canonical system (2.2) with three field rotation parameters,
$\lambda,$ $\beta,$ $\gamma,$ and suppose that (2.2) has three 
limit cycles around the origin. Then we get into some domain of the 
field rotation parameters being restricted by definite con\-di\-tions 
on two other parameters, $a$ and $c,$ corresponding to one of that 
cases of finite singularities which were considered in \cite{Gaiko}. 
We can fix both of these parameters putting, for example, $a=1$ 
and $c>1$ $(\nu=1)$ and supposing that system (2.2) has only
two finite singularities: a saddle and an anti-saddle. Thus, there 
is a domain bounded by two fold bifurcation surfaces forming 
a cusp bifurcation surface of multiplicity-three limit cycles in 
the space of the field rotation pa\-ra\-me\-ters $\lambda,$ 
$\beta,$ and $\gamma.$
    \par
The cor\-res\-pon\-ding maximal one-parameter family of
multiplicity-three limit cycles cannot be cyclic, otherwise there
will be at least one point cor\-res\-pon\-ding to the limit cycle
of multi\-pli\-ci\-ty four (or even higher) in the parameter
space. Extending the bifurcation curve of multi\-pli\-ci\-ty-four
limit cycles through this point and parameterizing the
corresponding maximal one-parameter family of
multi\-pli\-ci\-ty-four limit cycles by a field-rotation
para\-me\-ter, according to Theorem~4.2, we will obtain two
monotonic curves of multi\-pli\-ci\-ty-three and one, respectevely, 
which, by the Wintner--Perko termination principle (Theorem~4.1), 
terminate either at the origin or on a separatrix loop surrounding 
the origin. Since we know at least the cyclicity of the singular point 
which is equal to two (see \cite{BL}), we have got a contradiction 
with the termination principle stating that the multiplicity of limit cycles 
cannot be higher than the multi\-pli\-ci\-ty (cyclicity) of the singular point 
in which they terminate.
    \par
If the maximal one-parameter family of multiplicity-three limit
cycles is not cyclic, using the same principle (Theorem~4.2), this
again contradicts with the cyclicity of the origin \cite{BL} not admitting 
the multiplicity of limit cycles to be higher than two. This contradiction
completes the proof.\qquad $\Box$

\end{document}